\documentclass[letterpaper, 10 pt]{amsart}
\usepackage{geometry}                
\usepackage[parfill]{parskip}    
\usepackage{graphicx}
\usepackage{amssymb}
\usepackage{epstopdf}
\newcommand{\Rset}{\mathbb{R}}

\newtheorem{remark}{Remark}
\title{On Goal-Oriented Multiobjective Embedded Optimization for System Performance Assessment}
\author{Getachew K. Befekadu}
\keywords{Gradient-based optimization, nonlinear optimization, robust control}
\thanks{G. K. Befekadu is with the Department of Electrical \& Computer Engineering, 
Morgan State University, Baltimore, MD, USA.\\
E-mail: {\tt getachew.befekadu@morgan.edu}}%
\date{~~~June, 2010: Originally submitted as a technical report at the University of Notre Dame, Notre Dame, IN, USA}                                           

\begin{document}
\begin{abstract}
In this short note, we discuss a goal-oriented multiobjective optimization problem for system performance assessment. 
The objective function for such optimization problem, which is usually a composite of different performance indices corresponding 
to different operating conditions or scenarios in the system, is then posed as a goal-oriented multiobjective optimization problem. 
The (sub)-optimal solution(s) for such nonlinear optimization problem can be solved 
using a Sequential Quadratic Programming algorithm.  
\end{abstract} 
\maketitle

The dynamic performance of large system can often be enhanced by formulating the problem as a 
multiobjective dynamic embedded optimization problem (see e.g., \cite{1} - \cite{4}). In such optimization problem, 
the associated cost function, which usually embeds the evolving system dynamics as a feasible set constraint, 
can be used for assessing as well as improving the qualitative behavior of the system. 
For a large class of problems, it is sufficient to assume that the cost function is smooth, though the 
 underlying dynamic response is non-smooth due to intrinsic interactions between continuous dynamics 
 and discrete events in the system. A typical multiobjective dynamic embedded optimization can be 
 reformulated using the following objective function
\begin{align}
 \underset {x, \theta, t_f} {min} ~&J(x, \theta, t_f) \notag \\ 
 \quad  s.~t. \quad x(t) &= \Phi(t, \theta, x_0) \notag \\
  &\quad  x \in \Omega \subseteq \Rset^{n} \quad \textbf {(state space constraint)}\notag \\
 &\quad  \theta \in \Theta \subseteq \Rset^{p} \quad \textbf {(parameter space constraint)} \label{1}
\end{align}
where
\begin{align}
 J(x, \theta, t_f) = \varphi(x(t_f), \theta, t_f) + \int_{t_0}^{t_f} \psi(x(t), \theta, t) dt \label{2}
\end{align}
and $\theta$ are design parameters (i.e., controller parameters to be optimized in 
the system) and $\Phi(t, \theta, x_0)$ defines the evolution of the dynamic system. 
Moreover, $t_f$ is the final time and its adjustability is usually problem specific.

Though, the above problem formulation attempts to solve the optimal parameters for a particular
 operating condition or scenario, these controllers when implemented in the system do not guarantee 
 robustness for other operating conditions \cite{1}, \cite{2} and \cite{4}. In this short note, we discuss an optimization problem that 
 considers several operating conditions or system scenarios within the optimization framework of \eqref{1} 
 so as to guarantee robustness in the system for these envisaged operating conditions. 
 Thus, the optimization problem that considers several operating conditions can be reformulated 
 using the following goal-oriented composite objective function
\begin{align}
 \ {min}~&\boldsymbol{\gamma} \notag \\ 
 \quad  s.~t. \quad x(t) &= \Phi_i(t, \theta, x_0) \notag \\
 &\quad  x \in \Omega \subseteq \Rset^{n} \quad \textbf {(state space constraint)}\notag \\
 &\quad  \theta \in \Theta \subseteq \Rset^{p} \quad \textbf {(parameter space constraint)} \notag \\
 &J_i(x, \theta, t_f) - \gamma_i w_i \le J_{i}^*(x, \theta_{i}^*, t_f) \notag \\
 &w_i > 0, \notag \\
 &\gamma_i \quad \textbf {(are~unrestricted~scalar~variables)} \notag \\
 &i = 1, 2, \ldots, N \quad \textbf{(operating conditions)}  \label{3}
\end{align}
where
\begin{align}
 \ J_i(x, \theta, t_f) = \varphi_i(x(t_f), \theta, t_f) + \int_{t_0}^{t_f} \psi_i(x(t), \theta, t) dt  \label{4}
\end{align}
is the $i$th-objective function, $w_i$ is the weighting factor, $\Phi_i(t, \theta, x_0)$ defines the 
evolution of the dynamic system for the $i$th-operating condition, and the vector $\boldsymbol{\gamma}$ is 
given by $[\gamma_1, \gamma_2, \ldots, \gamma_N]^T$. Furthermore, the parameter vector $\theta_{i}^*$ and 
the corresponding achievable objective value $J_{i}^*$ are the (sub)-optimal solution for the $i$th-operating condition.

The problem formulation in (3) requires conceptually a two stage optimization strategy
 to solve the problem, i.e.,
 
\setlength {\hangindent}{20pt}
\noindent
\quad-~\textbf{First} solving the (sub)-optimal solution of the parameter vector $\theta_{i}^*$ and the corresponding achievable
 objective value $J_{i}^*$ (i.e., $i$th-objective goal) for the $i$th-operating 
 condition (assuming there are $N$ operating conditions), and
 
\setlength {\hangindent}{20pt}
\noindent
\quad-~\textbf{Second} solving a single parameter vector $\theta$ for these $N$ operating 
conditions using goal-oriented composite objective function of \eqref{3}.

The above two stage strategy, besides solving a single parameter vector $\theta$ for the corresponding
 operating conditions, gives quantitative information about the degree of under - and over - attainment of the goals $J_{i}^*$ 
 for different weighting coefficients $w_i$.  This further gives a qualitative measure about the robustness of 
 the controller(s) since the latter designing step (or controller(s) parameter optimization) accounts 
 different operating conditions and possibly different system scenarios simultaneously. 
 Moreover, in the above optimization problem, it is possible to incorporate structurally 
 different systems (i.e., systems with different topology structures corresponding to different operating conditions). 
However, the structure of the dynamic state-variable space for each system must remain topologically equivalent \cite{1}.

Therefore, the goal-oriented multiobjective optimization problem is based on the concept of gathering the values of all objectives 
(i.e., for each operating condition) into a single value (see references \cite{5} and \cite{6}). 
Then, the optimization problem can be solved using a multiobjective gradient-based optimization algorithm such as a Sequential 
Quadratic Programming (SQP) algorithm. Moreover, the weighting values in the composite objective 
function can be appropriately chosen to reflect the trade-off among the individual objectives. In contrast to other algorithms, 
the approach discussed in this note can provide solutions for a set of weighting values that represents or reflect $a~priori$ known system operations. 
Further analysis of these weighting values guides how to choose proper weightings for the individual objectives 
that could be used for the composite objective function.

\begin{remark}
If the structure of the dynamic state-variable varies during the course of its operation (i.e., the corresponding system is topological different from time to time), then the above problem formulation may not be directly applicable.
\end{remark}

\begin{remark}
The weighting factors are usually assigned positive values that sum to unity, i.e., $1=\sum_{i=1}^N{w_i}~and~w_i > 0$.
\end{remark}

\end{document}